\documentclass{article}
\usepackage{amssymb}
\usepackage{float}
\usepackage{color}
\usepackage{amsthm}
\usepackage{mhsetup}
\usepackage{mathtools}
\usepackage{tcolorbox}
\mathtoolsset{showonlyrefs} 
\numberwithin{equation}{section}

\newtheorem{Theorem}{Theorem}[section]

\theoremstyle{definition}

\theoremstyle{remark}
\newtheorem{Remark}[Theorem]{Remark}

\begin{document}
\title{Reparametrizing the relativistic Kepler equation: a bridge to Levi-Civita-type models\thanks{The authors are supported by INdAM-GNAMPA, Italy.}
}
\author{Alberto Boscaggin and Walter Dambrosio}
\date{}
\maketitle
\begin{abstract}
We establish a link between different relativistic variants of the Kepler problem. In particular, we show that solutions of the special relativistic model with fixed energy can be reparameterized as solutions of a generalized Kepler equation with an additional $1/r^2$ term in the gravitational potential. This yields a dynamics of the same type as the Levi-Civita correction, up to different coefficients.
\end{abstract}
\medskip



\section{Introduction}

In classical mechanics, the Kepler problem describes the motion of a body of mass $m$ in the gravitational field generated by an heavy body of mass 
$M$ fixed at the origin. Denoting by $G$ the gravitational constant and setting $\alpha = G M m$, the motion equation writes as
\begin{equation}\label{eq-kep}
m \ddot x = - \alpha \frac{x}{\vert x \vert^3}, \qquad x \in \mathbb{R}^3 \setminus \{0\}.
\end{equation} 
As is well known, bounded solutions of this equation describe
ellipses with one focus at the origin. In particular, the model predicts no
relativistic effects, such as the slow precession observed in the motion of
some planets.

A correct explanation of such effects is in fact provided by General Relativity.
More precisely, the motion of a test particle in the gravitational field of a
spherically symmetric body is described by timelike geodesics of the
Schwarzschild metric. Restricting to planar motion and exploiting the
conservation of energy and angular momentum, the radial dynamics can be
expressed in terms of an effective potential of the form
\begin{equation}\label{eq:VeffGR}
V_{\mathrm{eff}}(r)
=
-\frac{\alpha}{r}
+
\frac{L^2}{2mr^2}
-
\frac{GM L^2}{m c^2 r^3},
\end{equation}
where $L$ denotes the angular momentum of the particle and $c$ is the speed of light, cf. for instance \cite[Chapter~12.3~B]{GoPoSa02}.
The first two terms coincide with the classical Newtonian effective potential,
while the last one represents the leading relativistic correction. This term,
which vanishes in the non-relativistic limit $c \to \infty$, is proportional
to $r^{-3}$ and produces a small deviation from closed Keplerian ellipses.
In particular, it accounts for the well-known perihelion precession of
bounded orbits.

A different relativistic correction to the Newtonian model was proposed by
Levi-Civita \cite{LC-28}, starting from considerations inspired by geometrical approach of General Relativity. Indeed, Levi-Civita introduced a metric that describes the dynamics of a moving body subject to a gravitational potential, with an approximations of an order higher than $v^2/c^2$ (where $v$ is the velocity of the moving body). He proved that geodesic motions with respect to this metric can be approximated by Newtonian trajectories under the attraction a
modified effective potential of the form 
\begin{equation}\label{eq:VeffLC}
V_{\mathrm{eff}}(r)
=
-\frac{\alpha \left(1 + \frac{4  h}{mc^2}\right)}{r}
+
\frac{L^2}{2mr^2}
-
\frac{3 G^2 M^2 m}{c^2 r^2},
\end{equation}
where $h$ denotes the (approximated) energy of the particle.
In this expression, relativistic corrections appear both in the Newtonian
potential (which is multiplied by the factor $1 + 4h/mc^2$) and in the additional contribution
proportional to $r^{-2}$. Notice that, unlike the Schwarzschild
correction, the extra term is proportional to $r^{-2}$ rather than $r^{-3}$.

Both of these relativistic corrections emerge at the level of the
reduced radial equation obtained after restricting to planar motion
and exploiting the conservation of energy and angular momentum.
The corresponding effective potentials contain additional
inverse–power contributions, corresponding to a central force
problem of the type
\begin{equation}\label{eq:relkep}
m \ddot x = -\widehat\alpha \frac{x}{|x|^3} - \widehat\beta \frac{x}{|x|^{\ell}},
\end{equation}
with 
\[
\ell = 5, \qquad \widehat\alpha = \alpha, \qquad \widehat\beta = \frac{3GM L^2}{m c^2}
\]
in the case of the Schwarzschild correction and
\begin{equation} \label{eq-coeffLCoriginale}
\ell = 4, \qquad \widehat\alpha = \alpha \left( 1  +  \frac{4 h}{mc^2}\right), \qquad \widehat\beta = \frac{6 G^2 M^2 m}{c^2}
\end{equation}
for the Levi-Civita one. Note that this identification should be understood only at a formal
level. Indeed, the coefficients appearing in \eqref{eq:relkep} depend
on the first integrals of the motion and therefore cannot be interpreted
as fixed parameters of an underlying three–dimensional force law. Despite this, equations of the form \eqref{eq:relkep} with $\ell = 4$ or $\ell = 5$
are often referred to as relativistic Kepler problems and have been considered, among others, in \cite{AmBe90,LaLlNu91,NuCaLl91}.

On the other hand, more recently, a formulation of the Kepler problem has
also been proposed within the framework of Special Relativity \cite{AnBa71,Bo04,Ji13,LeMoPP,MuPa06}.
Precisely, keeping unchanged the classical Newtonian potential and replacing the classical
momentum by the relativistic one leads to the equation
\begin{equation}\label{eq-keprelintro}
\dfrac{\mathrm{d}}{\mathrm{d}t}\left(\dfrac{m\dot{x}}{\sqrt{1-\vert \dot{x}\vert ^2/c^2}}\right)
= -\alpha \, \dfrac{x}{\vert x\vert ^3}, \qquad x \in \mathbb{R}^3 \setminus \{0\}.
\end{equation}
This model has serious theoretical limitations, since it does not come from a fully relativistic theory of
gravitation. On the other hand, in contrast 
with the previous constructions which arise only at the
level of the reduced radial dynamics through effective potentials, it provides a genuine equation of motion in $\mathbb{R}^3$. 

The aim of this brief note is to point out a seemingly unnoticed simple relationship between equations \eqref{eq-keprelintro} and \eqref{eq:relkep}. More precisely, we prove that
\begin{quote}
\textit{solutions to equation \eqref{eq-keprelintro} with (relativistic) energy equal to $h$ can be reparameterized as solutions to an equation of the type \eqref{eq:relkep} with 
\begin{equation} \label{eq-coeffrel}
\ell = 4, \qquad \widehat\alpha = \alpha \left( 1  +  \frac{h}{mc^2}\right), \qquad \widehat\beta = \frac{ G^2 M^2 m}{c^2}.
\end{equation}}
\end{quote}
Note that this equation is exactly of the same type as the one proposed by Levi-Civita, up to different constant factors. As observed in the introduction \cite{LaLlNu91}, these constants are the ones which are obtained in \cite[pp. 217--218]{Be42} when studying the trajectories of \eqref{eq-keprelintro} using polar coordinates.


\section{Statement and proof of the result}\label{sec-2}

Let us consider the general relativistic problem
\begin{equation} \label{eq-rel}
\frac{\rm d}{{\rm d}t}\left(\frac{m\dot{x}}{\sqrt{1-|\dot{x}|^2/c^2}}\right)=\nabla V(x),\qquad x\in\Omega, 
\end{equation}
where $m, c > 0$, $\Omega \subset \mathbb{R}^n$ is an open set and $V:\Omega \subset \mathbb{R}^n\to \mathbb{R}$ is a function of class $C^1$.

\medskip

A simple computation proves that for every solution $x:I\subset \mathbb{R}\to \Omega$ of \eqref{eq-rel} there exists $h\in \mathbb{R}$ such that
\begin{equation}\label{eq-consenergy}
mc^2\left(\frac{1}{\sqrt{1-\dfrac{|\dot{x}(t)|^2}{c^2}}}-1\right)-V(x(t))=h,
\end{equation}
for every $t\in I$. The expression in the left-hand side of \eqref{eq-consenergy} is called the \textit{relativistic} energy of the solution.

Let us now define the set
\begin{equation} \label{eq-defomegah}
{\Omega}_h=\{x\in \Omega:\ V(x)+h \geq 0 \}.
\end{equation}
From \eqref{eq-consenergy} it is immediate to see that the necessary condition
\begin{equation} \label{eq-condnec}
V(x(t))+h\geq 0,\quad \forall \ t\in I, 
\end{equation}
holds true, thus implying that solutions are confined in the region $\Omega_h$.

\medskip

From now on, we will fix $h\in \mathbb{R}$. Let us consider the function $Z_h:\Omega \to \mathbb{R}$ defined by
\begin{equation} \label{def-Z}
Z_h(x)=V(x)+\dfrac{1}{2mc^2}\, (V(x)+h)^2,
\end{equation}
for every $x\in \Omega$. For future reference, let us observe that for every $x\in \Omega$ we have
\begin{equation} \label{eq-contoZ}
Z_h(x) + h = (V(x)+h)\, \left(1+\dfrac{1}{2mc^2}\, (V(x)+h)\right). 
\end{equation}

\medskip

Consider also the second order equation 
\begin{equation} \label{eq-duepunti}
m z''=\nabla Z_h(z)
\end{equation}
and trivially observe that for every solution $z$ of \eqref{eq-duepunti} there exists $\lambda \in \mathbb{R}$ such that
\begin{equation} \label{eq-energiaduepunti}
\dfrac{1}{2} m |z'|^2-Z_h(z)=\lambda.
\end{equation}

Our aim is to prove that there is a relation between solutions of \eqref{eq-rel} with relativistic energy $h$ and solutions of \eqref{eq-duepunti} with energy $h$. More precisely, the following results hold true.

\begin{Theorem}\label{teo-dareladuepunti}
	Let $x:[0,T]\to \Omega_h$ be a solution of \eqref{eq-rel} with relativistic energy $h$. Then, there exist $S>0$ and a change of variable $\chi: [0,S] \to [0,T]$ such that the function $z=x\circ \chi$ is a solution of \eqref{eq-duepunti} with energy $h$.
	
	\smallskip
	
	Conversely, let $z:[0,S]\to \Omega_h$ be a solution of \eqref{eq-duepunti} with energy $h$. Then, there exist $T>0$ and a change of variable $\eta: [0,T] \to [0,S]$ such that the function $x=z\circ \chi$ is a solution of \eqref{eq-rel} with relativistic energy $h$.
\end{Theorem}

\begin{Remark}
The choice of the Kepler potential $V(x) = \alpha / \vert x \vert$, with $\alpha = GMm$, provides
\begin{align*}
Z_h(x)  = \frac{\alpha}{\vert x \vert} + \frac{1}{2mc^2} \left(\frac{\alpha}{\vert x \vert} + h \right)^2 = GMm \left(1 + \frac{ h}{mc^2} \right)\frac{1}{\vert x \vert} + \frac{G^2 M^2 m}{2c^2} \frac{1}{\vert x \vert^2} + \frac{h^2}{2mc^2}
\end{align*}
and so the result mentioned in the Introduction follows. \hfill $\diamond$
\end{Remark}

\begin{Remark}\label{rem-duezone}
	Let us notice that the conservation of energy 
	\[
	\dfrac{1}{2}\, m\, |z'(s)|^2-Z_h(z(s))=h,\quad \forall \ s\in [0,S],
	\]
	for a solution $z$ of \eqref{eq-duepunti} implies the sign condition
	\[
	Z_h(z(s))+h\geq 0, \quad \forall \ s\in [0,S].
	\]
	Recalling \eqref{eq-contoZ}, we deduce that this condition is satisfied when
	\[
	V(z(s))+h+2mc^2\leq 0 \quad \mbox{or} \quad V(z(s))+h\geq 0,\quad \forall \ s\in [0,S].
	\]
	As a consequence, taking into account the continuity of $z$, we infer that solutions of \eqref{eq-duepunti} with energy $h$ are confined in $\Omega_h$ or in the region
	\[
	\Sigma_h=\{x\in \mathbb{R}^n:\ V(x)+h+2mc^2\leq 0\}.
	\]
	The crucial point in Theorem \ref{teo-dareladuepunti} is that the only solutions of \eqref{eq-duepunti} with energy $h$ which correspond to solutions of \eqref{eq-rel} \textit{with energy $h$} are the ones whose support lies in $\Omega_h$. In some sense, this is not unexpected, since a reparametrization does not modify the support of the solutions and, in order to have a solution with relativistic energy $h$, we necessarily need to start with a solution of \eqref{eq-duepunti} lying in $\Omega_h$.
	
	For completeness, let us remark that with a simple motification of the proof of Theorem \ref{teo-dareladuepunti} it is possible to show that solutions of \eqref{eq-duepunti} with energy $h$ and support in $\Sigma_h$ correspond to solutions of \eqref{eq-rel} with energy $h+2mc^2$, as expected. \hfill $\diamond$
\end{Remark}

Theorem \ref{teo-dareladuepunti} is basically a consequence of the relativistic Maupertuis principle recently proved in \cite{BDM}. However, for the reader's convenience we give here a more direct proof.

We start with the proof of the first statement. Given $x$ as in the statement, let us define
\begin{equation}
\zeta(t)=mc^2\, \int_0^t \dfrac{1}{V(x(\xi))+h+mc^2}\, d\xi,
\end{equation}
for every $t\in [0,T]$. 

Please notice that $\zeta$ is well-defined because $x(t)\in \Omega_h$, for every $t\in [0,T]$, implies that $V(x(t))+h+mc^2\geq mc^2>0$, for every $t\in [0,T]$. From this inequality we also deduce that $\zeta$ is strictly increasing (and invertible).

Let now define $S=\zeta(T)>0$ and let $\chi:[0,S]\to [0,T]$ be the inverse of $\zeta$. Let $z=x\circ \chi$.

\noindent
A simple computation shows that
\[
\begin{array}{l}
\displaystyle z'(s)=\dot{x}(\chi(s)) \, \dfrac{d\chi}{ds}=\dot{x}(\chi(s)) \, \dfrac{1}{\dot{\zeta} (\chi(s))}=\dot{x}(\chi(s)) \, \dfrac{V(x(\chi(s))+h+mc^2}{mc^2},\\
\end{array}
\]
for every $s\in [0,S]$.
We then deduce that
\begin{equation} \label{eq-contoenergia1}
\begin{array}{l}
\displaystyle
\dfrac{1}{2}\, m\, |z'(s)|^2 = \dfrac{1}{2m c^4}\, |\dot{x}(\chi(s))|^2 \, (V(x(\chi(s))+h+mc^2)^2,
\end{array}
\end{equation}
for every $s\in [0,S]$. On the other hand, from \eqref{eq-consenergy} and \eqref{eq-contoZ}, we infer that
\[
|\dot{x}(t)|^2=2mc^4\, \dfrac{(V(x(t))+h)\, \left(1+\dfrac{1}{2mc^2}\, (V(x)+h)\right)}{(V(x(t))+h+mc^2)^2}=2m c^4 \, \dfrac{Z_h(x(t))+h}{(V(x(t))+h+mc^2)^2},
\]
for every $t\in [0,T]$. By replacing in \eqref{eq-contoenergia1}, we obtain
\[
\dfrac{1}{2}\, m\, |z'(s)|^2 = Z_h(x(\chi(s)))+h=Z_h(z(s))+h,
\]
for every $s\in [0,S]$, thus proving the validity of the conservation law \eqref{eq-energiaduepunti}. From this relation, we immediately deduce \eqref{eq-duepunti} and this completes the proof of the first statement.

\medskip

As far as the proof of the second statement is concerned, assume that $z$ is a solution of \eqref{eq-duepunti} with energy $h$ and such that $z(s)\in \Omega_h$, for every $s\in [0,S]$. 

Consider the Cauchy problem
\[
\displaystyle \dot{\eta}=\dfrac{mc^2}{V(z(\eta))+h+mc^2},\quad \eta(0)=0
\]
and denote by $[0,\tau)$ the maximal interval of definition of the solution $\eta$. From the assumption $z(s)\in \Omega_h$, for every $s\in [0,S]$, we deduce that $0<\dot{\eta}(t)\leq 1,$ for every $t\in [0,\tau)$, thus implying that $\eta$ cannot blow-up in a finite time. From this fact, noticing also that the differential equation has no equilibria, we deduce that $\tau=+\infty$, that $\eta$ is strictly increasing on $[0,+\infty)$ and that $\eta$ goes to $+\infty$ as $t\to +\infty$.
As a consequence, there exists $T>0$ such that $\eta(T)=S$. Let now $x:[0,T]\to \mathbb{R}^n$ be defined by $x=z\circ \eta$. We then have
\[
\begin{array}{l}
\displaystyle \dot{x}(t)=z'(\eta(t))\, \dot{\eta}(t)=z'(\eta(t))\, \dfrac{mc^2}{V(z(\eta(t)))+h+mc^2},
\end{array}
\]
for every $t\in [0,T]$. From this relation, taking into account \eqref{eq-contoZ} and \eqref{eq-energiaduepunti}, we infer that
\begin{equation} \label{eq-contoenergia2}
\begin{array}{l}
\displaystyle
\dfrac{|\dot{x}(t)|^2}{c^2}= |z'(\eta(t))|^2 \, \dfrac{m^2c^4}{(V(z(\eta(t)))+h+mc^2)^2}=  \dfrac{2mc^2\, (Z_h(z(\eta(t)))+h)}{(V(z(\eta(t)))+h+mc^2)^2}\\
\\
\displaystyle =\dfrac{2mc^2\, (V(z(\eta(t)))+h)\, \left(1+\dfrac{1}{2mc^2}\, (V(z(\eta(t)))+h)\right)}{(V(z(\eta(t)))+h+mc^2)^2}=\dfrac{(V(x(t))+h)\, \left( V(x(t))+h+2mc^2\right)}{(V(x(t))+h+mc^2)^2}
\end{array}
\end{equation}
and
\begin{equation} \label{eq-contoenergia3}
\begin{array}{l}
\displaystyle
1-\dfrac{|\dot{x}(t)|^2}{c^2}=1-\dfrac{(V(x(t))+h)\, \left( V(x(t))+h+2mc^2\right)}{(V(x(t))+h+mc^2)^2} =\dfrac{m^2c^4}{(V(x(t))+h+mc^2)^2},
\end{array}
\end{equation}
for every $t\in [0,T]$. From this relation we plainly deduce that the conservation law \eqref{eq-consenergy} is satisfied. From this relation, we immediately deduce \eqref{eq-rel} and this completes the proof.

\noindent
A. Boscaggin\\
Dipartimento di Matematica ``Giuseppe Peano'', Università degli Studi di Torino\\
Via Carlo Alberto 10, 10123 Torino, Italy\\
alberto.boscaggin@unito.it
\medskip

\noindent
W. Dambrosio\\
Dipartimento di Matematica ``Giuseppe Peano'', Università degli Studi di Torino\\
Via Carlo Alberto 10, 10123 Torino, Italy\\
walter.dambrosio@unito.it

\end{document}